\documentclass[11pt]{amsart}
\usepackage{amssymb}
\usepackage{amsmath}
\usepackage{hyperref}
\usepackage{xcolor}

\newtheorem{theorem}{Theorem}[section]
\newtheorem{lemma}[theorem]{Lemma}

\newtheorem{cor}[theorem]{Corollary}
\newtheorem{conj}[theorem]{Conjecture}
\newtheorem{claim}[theorem]{Claim}
\newtheorem{construction}[theorem]{Construction}
\newtheorem{question}[theorem]{Question}

\newtheorem{definition}[theorem]{Definition}

\begin{document}
	\title[Treeings without matchings]{The measurable Hall theorem fails for treeings}
	
	\author{G\'abor Kun}
	
	

	\thanks{HUN-REN Alfr\'ed R\'enyi Institute of Mathematics + E\"otv\"os L\'or\'and University \\
		Budapest, Hungary \\
		Email: kungabor@renyi.hu\\
{\bf Acknowledgement.} The author thanks to the American Institute of Mathematics and the Fields Institute, where discussions have led to further applications. The author is particularly indebted to Ron Peled for suggesting the connection to couplings, and to L\'aszl\'o T\'oth and Zolt\'an Vidny\'anszky for their advices to improve this paper.
\\
{\bf Funding.}
This work has been supported by the Hungarian Academy of Sciences Momentum Grant no. 2022-58 and ERC Advanced Grant ERMiD}

	\begin{abstract}
		We construct, for every $d \geq 3$,  a $d$-regular acyclic measurably bipartite graphing that admits no measurable perfect matching, resolving a problem of Kechris and Marks and a number of further questions.
		
		A dense variant of our construction yields a coupling of two standard Borel probability measure spaces whose support contains no deterministic coupling, though the conditional probabilities of the coupling measure are atomless. This refutes a conjecture of Gurel-Gurevich and Peled.

	\end{abstract}
	
	\maketitle


	\section{Introduction}
	
	{\bf Measurable combinatorics} is concerned with graph theoretic problems arising in the setting of Borel graphs (see, e.g., Marks \cite{MICM}).
	The foundations of this theory have been laid by Kechris, Solecki and Todor\v{c}evi\'c \cite{KST}. In this area, the typical questions have the following form: given a Borel graph and a combinatorial problem (e.g., coloring, matching, graph decomposition) can we find a solution which is constructive in some sense, for example, it is measurable or Borel.
	
	In this paper we will focus on the measurable context. Recall that a {\it graphing} is a probability measure preserving (pmp) Borel graph, and a graphing is {\it measurably bipartite}, has a {\it measurable perfect matching}, etc., if it admits a  bipartition, perfect matching, etc., which are measurable w.r.t. the underlying probability measure, see Gaboriau \cite{GICM} or  Lov\'asz \cite{L}. Finally, essentially acyclic graphings are called {\it treeings}.

	One of the most important research directions in measurable combinatorics is the investigation of {\bf matchings}. Matchings have attracted a great deal of attention in the last decades, since equidecompositions can be phrased in terms of perfect matchings in Borel graphs, see M\'ath\'e \cite{Mathe} or Tomkowicz and Wagon \cite{TW}. Using such a rephrasing, the Banach-Tarski paradox shows that the existence of a (non-constructive) perfect matching in a Borel graph does not guarantee the existence of a measurable one \cite{BT} in general. Nevertheless, isolating situations in which nice matchings exist remains an extremely intriguing and fruitful direction of research.
	
	A natural subclass of Borel graphs to consider is the collection of \textbf{hyperfinite} ones. Indeed, measurably hyperfinite Borel graphs - and hence actions of amenable groups - do not admit paradoxical decompositions, which would easily rule out the existence of measurable perfect matchings. The study of such Borel graphs was also motivated by Tarski's famous circle squaring problem \cite{Tarski}, which has been resolved by Laczkovich \cite{Lsquaring} in the positive direction. Laczkovich' solution relies on finding a perfect matching in a concrete hyperfinite Borel graph, which arises from a $\mathbb{Z}^2$-action on $\mathbb{R}^2$ by translations. For this, he evokes the Axiom of Choice, hence the solution is not constructive. The problem of whether circle squaring was possible by nicer pieces remained open for decades until Grabowski, M\'ath\'e and Pikhurko found a measurable equidecomposition \cite{GMP}. Subsequently, Marks and Unger gave a Borel solution \cite{MU}, while recently M\'ath\'e, Noel and Pikhurko \cite{MNP} constructed a Jordan measurable equidecomposition.
	
	On the negative side, Laczkovich constructed a $2$-regular measurably bipartite graphing without a measurable perfect matching \cite{Lmatching}. Conley and Kechris \cite{CK} modified this to get a $d$-regular example of a graphing for every even $d$ (note that the latter are not acyclic). Recently, Bowen, Sabok and the author \cite{BKS} found
	a characterization of bipartite hyperfinite graphings admitting a measurable perfect matching in terms of ends. They showed that ``most" hyperfinite one-ended bipartite graphings with a perfect matching also admit a measurable perfect matching. The core obstruction is always a two-ended graphing (that is, of linear growth), as suggested by the known examples. Complementing the Conley-Kechris result, they also proved that every $d$-regular hyperfinite bipartite graphing admits a measurable perfect matching if $d$ is odd. 
	
	Now let us discuss the case of \textbf{non-hyperfinite} graphings. This class includes ergodic pmp actions of groups with Kazhdan's Property (T), and they play a crucial role in the papers of Margulis \cite{Margulis}, Sullivan \cite{S} and Drinfeld \cite{D} on the solution of the Banach-Ruziewicz problem on finitely additive invariant measures of the sphere.
	
	
 In their seminal paper, Lyons and Nazarov \cite{LN} showed that bipartite Cayley graphs of non-amenable groups admit a factor of iid perfect matching almost surely. Cs\'oka and Lippner \cite{CL} extended this to every Cayley graph of a non-amenable group. That is, with respect to the iid measure on the Bernoulli action of these groups, there is a measurable perfect matching.
	
 In contrast, in the Borel case Marks \cite{Marks} constructed $d$-regular acyclic Borel graphs without a Borel perfect matching for each $d>2$. This makes the following problem rather natural.
	
	\begin{question}
		\label{q:main} Let $d>2$. Does every $d$-regular acyclic graphing admit a measurable perfect matching?
	\end{question}

 The $3$-regular case of this question appears in the survey of Kechris and Marks (Problem 12.6, \cite{KM}), and another variant has been already mentioned by Klopotowski, Nadkarni, Sarbadhikari and Srivastava (2.7. Remark, \cite{KNSS}). The $3$-regular case was also one of the problems at the  Descriptive graph theory workshop at AIM \cite{AIM}.

 In this paper, we answer this question negatively. Our key novel tool is the usage of circulations and their characterization in terms of potentials (see Lov\'asz \cite{Lflow}, Lemma \ref{potential}). Recall that a {\it circulation} on a locally finite graph $G$ is a flow without sinks or sources. Observe that if $G$ is a bipartite $d$-regular graph and $M$ is a perfect matching, one can define a bounded non-zero circulation on $G$ as follows: let $A$ and $B$ be the partition classes, for $a \in A$, $b \in B, (a,b) \in E(G)$ set  $\varphi(a,b)=d-1$ if $(a,b) \in M$ and  $\varphi(a,b)=-1$ otherwise. Clearly, the preceding construction preserves measurability, so we can rule out the existence of measurable perfect matchings by disproving the existence of bounded measurable circulations. Indeed, our main result reads as follows.
	
	\begin{theorem} \label{main}
		For every $d \geq 2$ there exists a $d$-regular measurably bipartite treeing $T$
		such that  every absolutely integrable antisymmetric circulation on $T$ is zero almost everywhere. In particular, $T$ admits no measurable perfect matching.
	\end{theorem}
	
	Note the contrast with the case of factor of iid perfect matchings: as mentioned above, they almost surely exist on $d$-regular trees by Lyons and Nazarov \cite{LN}. Theorem~\ref{main} strengthens Marks' result \cite{Marks}: since $T$ does not admit a measurable perfect matching, it cannot admit a Borel one.
	
	Roughly speaking, we build our examples as inverse limits of finite graph sequences, diagonalizing against every circulation in the subsequent stage. This allows us a lot of flexibility and yields a number of further results. In what follows, we list the most important ones.
	
	\textbf{Couplings.}
	Our proof techniques turn out to be useful for denser structures, too.
	Couplings of standard Borel probability measure spaces\footnote{They had been named doubly stochastic measures, while recently they are often called permutons, though permutons are on $[0,1] \times [0,1]$, and, unlike in this paper, the linear ordering plays an important role in their study pioneered by Hoppen, Kohayakawa, Moreira, R\'ath, and Sampaio \cite{HKMRS}.} are classical objects of study. We may view a coupling measure as a fractional perfect matching. Recall that for finite bipartite graphs the existence of a fractional perfect matching, Hall's condition and the existence of a perfect matching are equivalent. Given two probability measure spaces $(X_1,\mathcal{A}_1,\lambda_1),(X_2,\mathcal{A}_2,\lambda_2)$, a relation
    $R \subseteq X_1 \times X_2$ and a coupling measure $\pi$ on $(X_1 \times X_2, \mathcal{A}_1 \otimes \mathcal{A}_2)$ with marginals $\pi^1=\lambda_1$ and $\pi^2=\lambda_2$, respectively, we say that the coupling $\pi$ {\it satisfies the relation} $R$ if $\pi(R)=1$.
    A {\it deterministic coupling} is a coupling satisfying the graph of a pmp mapping $S \to T$ (as a relation), i.e., a coupling such that the conditional probabilities a.s. consist of a single atom of measure one.
    Gurel-Gurevich and Peled made the following conjecture.
    \begin{conj}\cite{GP13}
    Given two standard Borel probability measure spaces $(X_1,\mathcal{A}_1,\lambda_1)$,
    $(X_2,\mathcal{A}_2,\lambda_2)$ and a Borel relation $R \subseteq X_1 \times X_2$, if there exists a conditionally atomless coupling satisfying $R$ then there exists a deterministic coupling satisfying $R$.
    \end{conj}

    They verified the conjecture for certain special cases related to Poisson thickening. Classical results of Bollob\'as and Varopoulos \cite{BV75} imply the positive direction, in the case when one of the measure spaces is countable. We refute the general conjecture.

	\begin{theorem} \label{coupling}
    There exists a standard Borel probability measure space $(X,\mathcal{A},\lambda)$,
    a Borel relation $R \subseteq X \times X$ and a conditionally atomless coupling $\pi$ on $(X \times X, \mathcal{A} \otimes \mathcal{A})$ satisfying $R$ such that there is no other coupling satisfying $R$. In particular, there is no deterministic coupling satisfying $R$.
	\end{theorem}

	Note that the coupling constructed must be an extreme point in the convex, closed set of couplings. Lindenstrauss \cite{L65} gave an interesting characterisation of the extreme points in the set of doubly stochastic measures, while Losert \cite{L82} has constructed an extreme doubly stochastic measure that is not supported on the union of countably many graphs of functions. Our example is a strengthening of the latter result.
	
	\textbf{Separating classes of local problems}. One can also use our main theorem to answer some questions arising from the recent study of the connection between descriptive combinatorics and {\it local algorithms} pioneered by Bernshteyn \cite{B}. Brandt, Chang, Greb\'{\i}k, Grunau, Rozho\v{n} and Vidny\'anszky (Question 1, \cite{BCGGRV}) asked about
	the relationship between classes of locally checkable labeling (LCL) problems admitting a factor of iid, and a measurable solution on $d$-regular trees, respectively. The above result shows that the LCL corresponding to perfect matchings separates these classes. In fact, we show that this separation is even possible by Schreier graphings of free pmp actions, see the remark at the end of this section.

	\textbf{Circulations in Schreier-graphings}. We construct a free pmp action of the free product of $d$-many two-element groups $(\mathbb{Z}/2\mathbb{Z})^{*d}$ without non-zero antisymmetric circulations.
	
	\begin{theorem}
		\label{action}
		For every $d \geq 3$ there exists a free pmp action of $(\mathbb{Z}/2\mathbb{Z})^{*d}$ such that the corresponding Schreier graphing $\mathcal{S}$ with the standard set of generators admits no measurable non-zero bounded antisymmetric circulation.
		
		Hence there is no free pmp $\mathbb{Z}$-action on any subset with positive measure whose Schreier graphing is a subgraphing of $\mathcal{S}$.
	\end{theorem}

    While $\mathcal{S}$ contains no Schreier graphing of a free pmp $\mathbb{Z}$-action as subgraphing, in contrast, T\'oth and the author observed that every graphing without finite components contains the Schreier graphing of a free $\mathbb{Z}$-action as a Lipschitz subgraphing, that is, the distance of adjacent vertices is bounded in the original graphing \cite{KT}.
	
	\textbf{Balanced orientations}. Recall that a {\it balanced orientation} is an orientation of the edges of a graph such that the in-degree is equal to the out-degree at every vertex.
	Thornton asked (Problem 3.4, \cite{T}) if regular graphings with even degree admit a measurable balanced orientation and proved it for regular graphings with expansion (Theorem 2.8, \cite{T}). Recently, Bowen, Sabok and the author \cite{BKS} proved this for hyperfinite, one-ended graphings.

	Conley \cite{C} provided a counterexample for $d=4$, while Bencs, Hru\v{s}kov\'a and T\'oth \cite{BHT} constructed for every $d$ a $2d$-regular, transitive graph that admits no factor of iid balanced orientation. Their example is two-ended and quasi-isometric to the bi-infinite line.
	
	Since a measurable balanced orientation of a $2d$-regular bipartite graphing naturally gives rise to a measurable circulation, Theorem \ref{action} yields the following corollary.
	
	\begin{cor} \label{balanced}
		For every even $d$ there exists a free pmp action of $(\mathbb{Z}/2\mathbb{Z})^{*d}$ such that the corresponding Schreier graphing with the standard set of generators admits no measurable balanced orientation. 	

        In particular, it is not the Schreier graphing of any free pmp action of the free group $\mathbb{F}_{d/2}$.
	\end{cor}

	This free action of $(\mathbb{Z}/2{\mathbb{Z}}\mathbb)^{*d}$ cannot even be isometric orbit equivivalent to any free pmp action of $\mathbb{F}_{d/2}$, see Joseph \cite{J} on isometric orbit equivalence.

    Let us remark that Thornton's result \cite{T} implies that the Schreier graphing of the Berno\-ulli shift of a non-amenable group always admits a measurable balanced orientation. Equivalently, the Cayley graph of the group admits a factor of iid balanced orientation. On the other hand, the construction in Corollary~\ref{balanced} is a free pmp action of a non-amenable group without a measurable balanced orientation. This shows that the class of LCL problems admitting a factor of iid solution strictly contains the class of LCL problems with a measurable solution, (partially) answering Question 1 in \cite{BCGGRV} in case of Schreier graphings, too.

    {\sc Organization of the paper.} In Section \ref{s:notation} we summarize the most important definitions and theorems we will use. In Sections \ref{s:inverse} and \ref{s:proofs} we give a brief outline and the proofs of our main results. Finally, in Section \ref{s:further} we list a number of further possible research directions.


	\noindent

	\section{Definitions, notation}
	
	\label{s:notation}
	
	Graph limit theory considers graphs whose vertex set is a standard Borel probability measure space $(J, \mathcal{A}, \lambda)$.
	Recall that {\it graphings}\footnote{While the names ``graphing" and ``treeing" were introduced by Adams \cite{A} referring to the representation of classes of  an equivalence relation
		by probability measure preserving (pmp) graphs, the usage of these expressions have shifted: we call pmp Borel graphs {\it graphings} following Gaboriau \cite{G} and Lov\'asz \cite{L}.}
	are locally finite Borel graphs whose edge set is the countable union of graphs of pmp partial functions on $J$, and essentially acyclic graphings are called {\it treeings}.
	The measure $\lambda$ induces a measure $\eta$ on the set of edges: given a Borel set of edges $S \subseteq E(G)$ define $\eta(S)=\frac{1}{2} \int_{V(G)} \deg^S(x) d \lambda(x)$ if it exists, where $\deg^S(x)$ denotes the number of neighbors of $x$ via $S$.
	
	In order to be able to handle both the sparse and dense cases, we will work in a more general framework introduced by Lov\'asz \cite{Lflow}. He called quadruples $(J, \mathcal{A}, \lambda, \eta)$, where
	the {\it node measure} $\lambda$ is a probability measure on the set of \emph{vertices}, which is the standard Borel space $(J,\mathcal{A})$, and the {\it edge measure} $\eta$ is a symmetric measure on $(J \times J, \mathcal{A} \otimes \mathcal{A})$, {\it double measure spaces}.
	In this paper we can also assume that the edge measure $\eta$ is $\sigma$-finite.

	Clearly, a graphing is a special case of this, where $\eta$ is the edge measure described above. Note that the measure $\eta$ admits an essentially unique disintegration, and for graphings the conditional measures will be supported on atoms of measure one, the neighbors of a vertex.
	
	In the case when $\lambda$ is the Lebesgue measure on $[0,1]$, $\eta$ is absolutely
	continuous with respect to $\lambda \times \lambda$, in symbols $\eta \ll \lambda \times \lambda$, and the Radon-Nikodym derivative takes values in $[0,1]$, we get the notion of \emph{graphons}.
	
    Recall that the \emph{marginals} of a measure $\eta$ on $J \times J$ are the measures $\eta^1$ and $\eta^2$ defined by $\eta^1(A)=\eta(A \times J)$ and $\eta^2(A)=\eta(J \times A)$ for all $A \in \mathcal{A}$. If both marginals of $\eta$ are equal to $\lambda$, that is, $\eta^1=\eta^2=\lambda$, then we call it a {\it Markov space}: in this case $(J^2, \mathcal{A}^2, \eta)$ is a \emph{coupling} of $(J,\mathcal{A},\lambda)$ with itself on $(J \times J, \mathcal{A} \otimes \mathcal{A})$. Markov spaces will be our main objects of study.
	
	

	A {\it circulation} on a graph $G$ is a map $\phi: E(G) \rightarrow \mathbb{R}$ that satisfies the \emph{flow equation} for every $x \in V(G)$, that is, for every $x$ the equality \[\sum_{y: (x,y) \in E(G)} \phi(x,y) = \sum_{y: (x,y) \in E(G)} \phi(y,x)\] holds.
    The definition of circulation extends to Markov spaces if $\phi$ is absolutely integrable with respect to the edge measure $\eta$, in this case the sums in the equality above are integrals of a.e. well-defined functions.
    We will call circulations {\it circulation functions} in order to distinguish them from
    the more general circulation measures to be introduced later.

	Every {\it signed measure} $\alpha$ is the difference of two finite measures with disjoint support: $\alpha=\alpha^+-\alpha^-$. The {\it variation}
	of $\alpha$ is the finite measure $|\alpha|=\alpha^++\alpha^-$, while the {\it total variation} denoted by $\|\alpha\|$ is the $|\alpha|$-measure of the base set.
    A signed measure on $J \times J$ is {\it symmetric} if $\alpha(A \times B) = \alpha(B \times A)$ for every $A, B \in \mathcal{A}$, and {\it antisymmetric} if $\alpha(A \times B) = -\alpha(B \times A)$ for every $A, B \in \mathcal{A}$.
    Every signed measure on $J \times J$ is the sum of a symmetric and an antisymmetric signed measure: $\alpha= \frac{\alpha+\alpha^T}{2}+\frac{\alpha-\alpha^T}{2}$.
	
    We will use the general concept of circulations over a $\sigma$-algebra $(J, \mathcal{A})$ still following Lov\'asz \cite{Lflow}. A signed measure $\alpha$ on the product $\sigma$-algebra $(J \times J, \mathcal{A} \otimes \mathcal{A})$ is called a {\it circulation measure} if it has equal marginals: $\alpha^1=\alpha^2$. This is equivalent to $\alpha(J \times S)=\alpha(S \times J)$ for every $S \in  \mathcal{A}$.
    Every symmetric signed measure on $J \times J$ is a circulation measure, and every circulation measure is the sum of a symmetric and an antisymmetric circulation measure.
    Observe that this generalizes absolutely integrable circulation functions on double measure spaces: to every absolutely integrable circulation function $\varphi: E(G) \rightarrow \mathbb{R}$ we can associate the circulation measure $\alpha$ defined by $\alpha(S) = \int_S \varphi(x,y) \ d\eta(x,y)$ for every $S \in \mathcal{A}$. The signed measure $\alpha$ given by this formula is a circulation measure if the function $\varphi$ is a circulation function.
	
    Note that, in general, the definition of a circulation measure on a double measure space does not depend on the edge measure.
    In the sparse case, for graphings, it would be sufficient for us to use circulation functions only.

    Our key observation is that for regular bipartite graphs the existence of a perfect matching implies the existence of a non-zero circulation function. This generalizes to Markov spaces and circulation measures.

    \begin{lemma} \label{observe}
    Let $\mathcal{M}=(J,\mathcal{A},\lambda,\eta)$ be a Markov space, $J=A \cup^*B$ a Borel partition, $R \subseteq A \times B \cup B \times A$ a symmetric Borel relation such that
    $\eta(R)=1$.
    \begin{enumerate}
    \item
    If there is a coupling $\alpha$ on $(J \times J, \mathcal{A} \otimes \mathcal{A})$ such that $\alpha(R)=1$ and $\alpha \neq \eta$, then there is a non-zero antisymmetric circulation measure $\beta$ such that $|\beta|(J \times J \setminus R)=0$.

    \item
    If $\mathcal{M}$ admits no non-zero antisymmetric circulation satisfying $R$ and $\eta$ is  not a deterministic coupling then there is no deterministic coupling satisfying $R$.
    \end{enumerate}
    \end{lemma}

    \begin{proof}
    $(1)$ First, we replace $\alpha$ by a symmetric $\alpha'$.
    We may assume that $\alpha|_{A \times B} \neq \eta|_{A \times B}$. Define $\alpha'$ by $\alpha'(S)=\alpha(S \cap A \times B)+\alpha(S^T \cap A \times B)$ for every Borel $S \subseteq J \times J$. Note that $\alpha'$ is symmetric and $\alpha' \neq \eta$.

    Next, define the circulation measure $\beta$ for every Borel $S \subseteq A \times B$ by

    $$\beta(S)=\begin{cases}
     \alpha'(S)-\eta(S) \text{ if }
     S \subseteq B \times A
    \\
    \eta(S)-\alpha'(S) \text{ if }
    S \subseteq A \times B
    \\
    0 \text{ if }
    S \subseteq A \times A \cup B \times B,
    \end{cases}$$

    this extends to an antisymmetric circulation measure. $\alpha' \neq \eta$, hence $\beta \neq 0$.

    $(2)$ If there was a deterministic coupling $\alpha$ satisfying $R$ then $\alpha \neq \eta$,
    and this would contradict $(1)$.
    \end{proof}

    We will use the following characterization of circulation measures by Lov\'asz (Lemma 4.2, \cite{Lflow}) in terms of \emph{potentials}. Recall that a measurable function $F:J \times J \to \mathbb{R}$ is a \emph{potential}, if
    there is a measurable function $f:J \to \mathbb{R}$ such that $F(x,y)=f(x)-f(y)$ for all $x,y \in J$.
	
	\begin{lemma} \label{potential}
		A signed measure $\alpha$ in a double measure space is a circulation if and only if for every bounded\footnote{The boundedness condition has been implicitly assumed but omitted in \cite{Lflow}.} measurable potential $F$ the equality \[\int_{J \times J} F(x,y) \ d \alpha(x, y) = 0\] holds.
	\end{lemma}
	For the sake of completeness we copy the proof from \cite{Lflow}.
	
	\begin{proof}
	The "if" part follows by applying the condition to the potential $\chi_A(x)-\chi_A(y):$ for any $A \in \mathcal{A}$ we have\\
    \[\alpha(A \times J)-\alpha(J \times A)=\int_{J \times J} \chi_A(x)-\chi_A(y) \ d \alpha(x,y)=0.\]

	To prove the converse, let $\alpha$ be a circulation measure, then for every potential $F(x,y)=f(x)-f(y)$, we have\\
	\[\int_{J \times J} f(x)-f(y) \ d \alpha(x,y) = \int_J f(x) \ d\alpha^2(x) - \int_J f(y) \ d\alpha^1(y)=0.\]
	\end{proof}

    Given two graphs $G$ and $H$ a {\it graph homomorphism} from $G$ to $H$ is a mapping $f: V(G) \to V(H)$ such that $(f(x),f(y)) \in E(H)$ for every edge $(x,y) \in E(G)$.
    The {\it tensor product} (also called categorical or direct product) of $G$ and $H$ is the graph with vertex set $V(G) \times V(H)$ and edge set $\{ ((g,h)(g',h')): (g, g') \in E(G), (h,h') \in E(H) \}$.

    For an undirected graph $G$ we denote by $Ori(G)$ the set of the orientations of the edge set of $G$. We think of an orientation as an antisymmetric map $E(G) \to \{-1,1\}$. Note that $|Ori(G)|=2^{|E(G)|}$.
	
	The set of the first $n$ positive integers will be denoted by $[n]=\{1, \dots ,n\}$.

 	\section{Inverse limits}

	The main ideas of the proof of Theorems \ref{main} and \ref{coupling} are essentially the same.
	We construct an inverse limit of a sequence of finite graphs $\{ G_n \}_{n=1}^{\infty}$. It is not hard to guarantee the basic properties of the inverse limit: if the sequence consists of regular graphs then the limit will be a Markov space, and if the sequence is bipartite then the limit will have a Borel bipartition (see Lemma \ref{basic} below).
	
	The question is how to ensure that the limit has only symmetric circulation functions (measures in the proof of Theorem \ref{coupling}). The key idea is utilizing Lemma~\ref{potential}:
	this shows that if for every $n$ and every orientation of $G_n$ there is a potential on $G_{n+1}$ such that for most pairs of adjacent vertices in $G_{n+1}$ the difference of the potential agrees with the value of the $G_n$-orientation of the image of the edge connecting them, then the limit object admits only symmetric circulation functions (measures). This allows a diagonalization argument for a construction without a non-zero antisymmetric circulation function (measure), see Lemma \ref{limit}.
		
	\label{s:inverse}
	First we state the lemma on the construction of inverse limits of graph sequences. Let $\{ G_n \}_{n=1}^{\infty}$ be a sequence of finite graphs and $f_n:V(G_{n+1}) \rightarrow V(G_n)$ be graph homomorphisms. We define the inverse limit of the sequence $\{ G_n ,f_n\}_{n=1}^{\infty}$ to be the graph $G$ with vertex set
	\[ V(G)=J=\big\{ (x_n)_{n=1}^{\infty}:\forall n \text{ } x_n \in V(G_n), f_n (x_{n+1})=x_n  \big\}\]
	and edge set
	\[E(G)=\{ (x, y): x,y \in J, \forall n  \ (x_n, y_n) \in E(G_n) \}.\]
	We endow $J$ with the topology inherited from the product topology of the discrete topologies on $V(G_n)$.
	
	By slightly abusing the notation if $f:V(G') \to V(G)$ is a mapping and $(u,v) \in E(G)$, we will denote by $f^{-1}(u,v)$ the set \[\{(u',v'): f(u')=u \text{ and } f(v')=v\}.\]
	
	A sequence $\{ G_n,f_n \}_{n=1}^{\infty}$ will be called \emph{proper}, if for every $n$
	\begin{enumerate}
		\item $|f_n^{-1}(u)|=\frac{|V(G_{n+1}|}{|V(G_n)|}$ for every $u \in V(G_n)$,
		\item $|f^{-1}_n(u,v) \cap E(G_{n+1})|=\frac{|E(G_{n+1})|}{|E(G_n)|}$ for every $(u,v) \in E(G_n)$,
		\item the graph $G_n$ is regular with degree at least two, and
		\item the graph $G_n$ is bipartite.
	\end{enumerate}
	
	Let $\mathcal{A}$ be the Borel $\sigma$-algebra on $J$, then $\mathcal{A}$ is generated by cylinders. Set
	\[\lambda\big(\{ x:
	x_n \in S \}\big):=\frac{|S|}{|V(G_n)|},\]
	for every $n$ and $S \subseteq V(G_n)$, this uniquely extends to $\mathcal{A}$.

    Similarly, consider the product $\sigma$-algebra $\mathcal{A}^2$ on $J^2$. Define $\eta$ by \\ \[\eta\Big(\{ (x, y): (x_n, y_n) \in Q \}\Big):=\frac{|Q \cap E(G_n)|}{|E(G_n)|},\] for every $n$ and $Q \subseteq V(G_n)^2$, this uniquely extends to
	$\mathcal{A}^2$.

    Finally, in order to disintegrate $\eta$, given $x \in J$ define the probability measure $\eta_x$ on $(J,\mathcal{A})$ on the cylinder sets by
    $\eta_x(S):=\frac{deg_S(x_n)}{deg(x_n)}$ 	for every $n$ and $S \subseteq V(G_n)$, this uniquely extends to $\mathcal{A}$.

    \begin{lemma} \label{basic}
		Let $\{ G_n,f_n \}_{n=1}^{\infty}$ be a proper sequence, and define $G, \lambda$ and $\eta$ as above.
		Then
		\begin{enumerate}
			
			\item $V(G)$ is a Polish space, $E(G)$ is a closed subset of $V(G) \times V(G)$.
			
			\item $(J, \mathcal{A}, \lambda)$ is a probability measure space.

			\item $(J^2, \mathcal{A}^2, \eta)$ is a probability measure space and $\eta(J^2 \setminus E(G))=0$.
			
			\item $\mathcal{D}=(J,\mathcal{A},\lambda,\eta)$ is a double measure space.
			
			\item
			$\eta^1=\eta^2=\lambda$, in other words, $\mathcal{D}$ is a Markov space.
			
			\item
			If for some $d$ every $G_n$ is $d$-regular then 
			  $(J,\mathcal{A},\lambda,d \eta)$ is a $d$-regular graphing.
			
			\item
			$G$ admits a clopen bipartition.

             \item $\{ \eta_x \}_{x \in J}$ is a disintegration of $\eta$, that is, for every Borel $S \subseteq J \times J$ we have \[ \eta(S) = \int_J \eta_x(\{y: (x,y) \in S\}) d\lambda(x).\]

		\end{enumerate}
	\end{lemma}
	
	\begin{proof}
		$(1)$ and $(2)$ are immediate.
		
		For $(3)$, first note that \[\eta(J^2 \setminus E(G))=\eta(\{ (x, y): \exists n \ (x_n, y_n) \notin E(G_n) \})= \]
		\[\eta(\bigcup_{n=1}^{\infty} \{ (x, y): (x_n, y_n) \notin E(G_n) \}) \leq
		\sum_{n=1}^{\infty} \eta(\{ (x, y): (x_n, y_n) \notin E(G_n) \})=0.\]

		$(4)$ The quadruple $(J,\mathcal{A},\lambda,\eta)$ is a double measure space, since $\eta$ is symmetric on the cylinder sets and hence on all sets in $\mathcal{A}^2$.
		
		$(5)$ Since every $G_n$ is regular $\eta(S \times J)=\eta(J \times S)=\lambda(S)$ holds for every cylinder set $S$, and hence for every $S \in \mathcal{A}$.
		
		$(6)$ Observe that for any measurable subset $S \subseteq J^2$ the equality \[\eta(S)=\frac{1}{d} \int_J deg^S(x) d \lambda(x)\] holds for every cylinder set $S$ and hence for every $S \in \mathcal{A}^2$. Since the measure $\eta$ is symmetric $(J,\mathcal{A},\lambda,d \eta)$ is a $d$-regular graphing.
		
		$(7)$ The bipartition of $G_1$ induces a clopen bipartition of $\mathcal{D}$.

        $(8)$ holds for every cylinder set $S \subseteq J \times J$ by the definition of $\eta_x$, and hence for every Borel set in $\mathcal{A} \times \mathcal{A}$.
	\end{proof}
	
    The next lemma encompasses the key observation. In the case of finite graphs, to a given non-zero circulation, one can associate an orientation $\mathcal{O}$ of the edges, based on the sign of the circulation. Clearly, there can be no potential $p$ that matches the value of the orientation on each edge (that is, $p(x)-p(y)=\mathcal{O}(x,y)$): this would contradict Lemma \ref{potential}.

    Now, when we have an inverse limit of finite graphs, we eliminate non-trivial circulation measures by approximating them on some finite stage and by adding potentials in the next stage which match the value of the corresponding orientation. Let us formalize this argument.

    \begin{definition}
    Given a proper sequence $\{ G_n,f_n \}_{n=1}^{\infty}$, a positive integer $n$, an orientation $\mathcal{O} \in Ori(G_n)$ and a mapping $p: V(G_{n+1}) \to \mathbb{R}$ we say that an edge $(u,v) \in E(G_{n+1})$ is $(p,\mathcal{O})-${\it bad} if $p(u)-p(v) \neq \mathcal{O}(f_n(u),f_n(v))$.
    A vertex of $V(G_{n+1})$ is $(p,\mathcal{O})$-{\it bad} if it is incident to a bad edge.
    The set of $(p,\mathcal{O})$-bad vertices is denoted by $B_{p,{\mathcal{O}}}$.
    \end{definition}

	\begin{lemma} \label{limit}
		Let $\{ G_n,f_n \}_{n=1}^{\infty}$ be a proper sequence and let $\mathcal{D}=(J,\mathcal{A},\lambda,\eta)$ be the inverse limit Markov space associated to it. Assume that for every $\varepsilon>0$ there is an $N$ such that for every $n>N$ and orientation $\mathcal{O} \in Ori(G_n)$ there exists a mapping $p_{\mathcal{O}}: V(G_{n+1}) \rightarrow \mathbb{R}$ such that
		\begin{enumerate}
			\item for every $(u,v) \in E(G_{n+1})$ we have $|p_{\mathcal{O}}(u)-p_{\mathcal{O}}(v)| \leq 1$,
		      \item \label{c:important}
            $|B_{p_{\mathcal{O}},\mathcal{O}}| \leq \varepsilon |V(G_{n+1})|$ holds for the set of $(p_{\mathcal{O}},\mathcal{O})$-bad vertices.
        \end{enumerate}
		
		Then for every antisymmetric circulation measure $\alpha$ on $\mathcal{A}^2$ if $|\alpha|^1 \ll \lambda, |\alpha|^2 \ll \lambda$ and $|\alpha|(J \times J \setminus E(G))=0$ then $\alpha=0$.
	In particular, the only coupling on $J \times J$ that satisfies $E(G)$ is $\eta$, so there is no such deterministic coupling, i.e., $E(G)$ does not contain a.e. a pmp perfect matching.
		
	\end{lemma}

	Condition (2) above can be made simpler if the graphs in the sequence have uniform degree. In this case it would be sufficient to bound, by $\varepsilon |E(G_{n+1})|$, the number of bad edges instead of the number of bad vertices. If the sequence does not have bounded degree then the limit is denser, hence the graph of a pmp mapping is an $\eta$-nullset,
	and we need a different condition that is sensitive to circulation measures of smaller support.
	
	\begin{proof}
		
		Towards a contradiction, we suppose that there exists such an antisymmetric $\alpha \neq 0$. Since $|\alpha|^1 \ll \lambda, |\alpha|^2 \ll \lambda$, the Radon-Nikodym theorem gives an $\varepsilon>0$ such that for any measurable subset $S \subseteq J$ if $\lambda(S) \leq \varepsilon$ then
		$|\alpha|(S \times J \cup J \times S)<\frac{\|\alpha\|}{3}$.
		
		Since $\alpha$ is a signed measure, by the Hahn decomposition theorem, there exists a measurable map $\mathcal{O}':J^2 \to \{-1,1\}$ such that
        \[\int_{J^2} \mathcal{O}'(x,y) \ d \alpha(x,y)=\|\alpha\|.\]

        Note that $\mathcal{O'}$ is antisymmetric on an $|\alpha|$-conull set, so we may assume that it is antisymmetric everywhere.
		Since $|\alpha|(J \times J \setminus E(G))=0$ we can choose an integer $n$ large enough such that there exists an orientation $\mathcal{O}$ of $E(G_n)$ that satisfies
		\[\int_{J^2} \mathcal{O}(x_n,y_n) \ d \alpha(x,y)> \frac{2\|\alpha\|}{3}.\]
		There is a corresponding $p_{\mathcal{O}}:V(G_{n+1}) \to \mathbb{R}$ satisfying the conditions of the lemma for $\varepsilon$. Consider the set $S$ of vertices with a $(p_{\mathcal{O}},\mathcal{O})$-bad $(n+1)$st coordinate: $S=\{x: x_{n+1} \in B_{p_{\mathcal{O}},\mathcal{O}}\}$.
		The set $S$ is measurable, and by \eqref{c:important} and the definition of $\lambda$, we have $\lambda(S) \leq \varepsilon$.
		
	On the one hand, \[\int_{J^2} p_{\mathcal{O}}(x_{n+1})-p_{\mathcal{O}}(y_{n+1}) \ d\alpha(x,y) = 0\] by Lemma \ref{potential}, since $(x,y) \mapsto p_{\mathcal{O}}(x_{n+1})-p_{\mathcal{O}}(y_{n+1})$ is a potential. 
		On the other hand,
		$$\int_{J^2} p_{\mathcal{O}}(x_{n+1})-p_{\mathcal{O}}(y_{n+1}) \ d \alpha(x,y) \geq $$
		$$\int_{J^2} \mathcal{O}(x_n,y_n) \ d \alpha(x,y) - \int_{S \times J \cup J \times S} 2 \ d |\alpha|(x,y).$$
		
		This leads to a contradiction, since $\lambda(S) \leq \varepsilon$, and hence
		$$\int_{J^2} \mathcal{O}(x_n,y_n) \ d \alpha(x,y) - \int_{S \times J \cup J \times S} 2  \ d|\alpha|(x,y)> \frac{2\|\alpha\|}{3}-2|\alpha|(S)=0.$$

    The second, "in particular" part follows from Lemma \ref{observe}: in order to apply $(2)$ we need that the measure $\eta$ is not a deterministic coupling on $J \times J$, and this follows from $(8)$ of Lemma \ref{basic}, since the degree of $G_n$ is at least two for every $n$.
	\end{proof}
	
	\section{The proofs}
	
	\label{s:proofs}

    In order to assure that the assumption of Lemma~\ref{limit} holds, when $G_n$ is given, we will consider several disjoint copies of a big subgraph (Construction~\ref{construction}) of a topological cover of $G_n$ (as a simplicial complex, where the covering corresponds to the homomorphisms of its fundamental group to $\mathbb{Z}^{Ori(G_n)}$ associated to the orientations of $G_n$) with small boundary, and we add a small number of vertices and edges to make it regular and bipartite. A coordinate of the cover gives us the corresponding potential on $G_{n+1}$, and the error, where the difference of the potential on an edge does not agree with the orientation of the image in $G_n$ can occur only on the additional edges.
 	
	\begin{construction} \label{construction}
		Let $G$ be a finite, undirected, $d$-regular graph and $N \in \mathbb{N}$.
		Define the graph $F=F(G,N)$ as follows.
		$$\begin{displaystyle} V(F)= V(G) \times \Pi_{\mathcal{O} \in Ori(G) } [N] \end{displaystyle}.$$
		Let $f: V(F) \rightarrow V(G)$ and
		$p_{\mathcal{O}}: V(F) \rightarrow [N]$ for
        $\mathcal{O} \in Ori(G)$ denote the projections.
		Set
		\[E(F)=\{ (x,y): x,y \in V(F), (f(x),f(y)) \in E(G), \]
	    \[ \forall \mathcal{O} \in Ori(G) \text{ }p_{\mathcal{O}}(x)-p_{\mathcal{O}}(y)=\mathcal{O}(f(x),f(y))\}.\]
	\end{construction}
	
	We state the basic properties of this construction.
	
	\begin{lemma} \label{nextaction}
		Let $G$, $N$, $F=F(G,N)$, $f$ and $p_{\mathcal{O}}$ be as above.
		
		\begin{enumerate}
			
			\item
			The mapping $f$ is a graph homomorphism.
			
			\item
			Every $u \in V(F)$ has degree at most $d$, and if $p_{\mathcal{O}}(u) \notin \{1,N\}$ for every orientation $\mathcal{O} $ then the degree of $u$ is $d$.
			In particular, there are at most $\frac{2|Ori(G)|}{N}|V(F)|$ vertices with degree strictly less than $d$.
			
			\item
			$|f^{-1}(u)|=N^{|Ori(G)|}$ for every $u \in V(G)$.
			
			\item
			$|f^{-1}(u,v) \cap E(F)|= (N-1)^{|Ori(G)|}$ for every $(u,v) \in E(G)$.
			
			\item
			$|\{s \in f^{-1}(u) : \nexists t  \in f^{-1}(v), (s,t) \in E(F) \}|=  N^{|Ori(G)|} - (N-1)^{|Ori(G)|}$
			for every $(u,v) \in E(G)$.
		\end{enumerate}
	\end{lemma}
	
	\begin{proof}
		$(1), (2)$ and $(3)$ are straightforward consequences of the construction.
		In order to prove $(4)$ note that, given an edge $(u,v) \in E(G)$, if $(s,t) \in f^{-1}(u,v) \cap E(F)$ every orientation $\mathcal{O}$ determines whether $p_\mathcal{O}(s)-p_{\mathcal{O}}(t)$ equals to $1$ or $-1$.
		This gives $(N-1)$ possible values of the pair $(p_{\mathcal{O}}(s),p_{\mathcal{O}}(t))$ for every orientation $\mathcal{O}$.

        Note that the quantity in $(3)$ is the sum of those in $(4)$ and $(5)$: for every $s \in f^{-1}(u)$ there is at most one $t$ such that $(s,t) \in f^{-1}(u,v)$. Hence $|f^{-1}(u)|=|f^{-1}(u,v) \cap E(F)|+|\{s \in f^{-1}(u) : \nexists t  \in f^{-1}(v), (s,t) \in E(F) \}|$.
	\end{proof}
	
	Now we have all the tools to prove the theorems.
	\begin{proof}[Proof of Theorem~\ref{main}]
		We will define a sequence of finite $d$-regular bipartite graphs $\{G_n\}_{n=1}^{\infty}$ recursively.
		Let $G_1=K_{d,d}$ be the complete bipartite graph on $d+d$ vertices.
		When $G_n$ is defined, let $F_n=F(G_n,N_n)$, where $N_n= |Ori(G_n)|2^{n+2}$ and $f:F_n \to G_n$ as in Construction \ref{construction} above. Note that the graph $F_n$ is not regular. Below we will slightly modify it to make sure that this holds as well.
		
		Take $d^2$ vertex-disjoint isomorphic copies of $F_n$ and for every $d^2$-tuple belonging to a vertex $x \in V(F_n)$ add
		$(d-deg_{F_n}(x))d$ extra vertices each adjacent to $d$ distinct vertices from the $d^2$-tuple in order to get a $d$-regular graph. We will refer to the vertices in the copies of $F_n$ as {\it internal} vertices, and the rest of the vertices
        as {\it external} vertices.
		Formally, \[V(G_{n+1})=  [d^2] \times V(F_n) \cup \{ (x,j): x \in V(F_n), j \in [d(d - deg_{F_n}(x))] \},\] \\
		and \[E(G_{n+1})=\{((i,x),(i,y)): (x,y) \in E(F_n), i \in [d^2] \} \cup \] \[\big \{ \{(i,x),(x,j)\}: \]\[x \in V(F_n), \exists k \in [d-deg_{F_n}(x)], \ell \in [d],
		kd^2-i = jd-\ell\}.\]

		We extend the homomorphism $f: F_n \to G_n$ from the copies of $F_n$ to $G_{n+1}$ to a map $f_n$, by mapping $(x,j)$
		to the $m$th neighbor of $x$, where $j \equiv m \mod d$. The pre-image of every vertex of $G_n$ has the same size: this holds for the intersection with internal vertices by $(3)$ of Lemma~\ref{nextaction}, and for the intersection with external vertices by $(5)$ of Lemma~\ref{nextaction} and the definition of $f_n$.

        And $(4)$ of Lemma~\ref{nextaction} also implies that the intersection of the pre-image of every edge in $G_n$ with the set of edges connecting internal vertices has the same size, while $(5)$ shows it for the set of edges incident to an external vertex. Hence we have constructed a proper sequence, the assumptions of Lemma~\ref{basic} hold, and the limit is a $d$-regular graphing.

        For an orientation $\mathcal{O} \in Ori(G_n)$ let $\overline{p}_{\mathcal{O}}$ be the map on
        $V(G_{n+1})$ defined by
		\[\overline{p}_{\mathcal{O}}(i,x)=\overline{p}_{\mathcal{O}}(x,j)=p_{\mathcal{O}}(x).\]
     	In order to apply Lemma~\ref{limit}
		we estimate the number of $(\overline{p}_{\mathcal{O}},\mathcal{O})$-bad vertices.
        Note that $(\overline{p}_{\mathcal{O}},\mathcal{O})$-bad vertices are exactly the external vertices and internal vertices with an external neighbor, hence the set of $(\overline{p}_{\mathcal{O}},\mathcal{O})$-bad vertices is the same for every $\mathcal{O} \in Ori(G_n)$. The number of internal bad vertices is less than $2 \cdot \frac{|Ori(G_n)|}{N_n}|V(G_{n+1})|$ by $(2)$ of Lemma~\ref{nextaction}. All edges from external vertices go to bad internal vertices, so the external ones contribute at most half of the bad vertices by regularity. We conclude that
        \[|B_{\overline{p}_{\mathcal{O}},\mathcal{O}}| \leq \frac{4|Ori(G_n)|}{N_n} |V(G_{n+1})|=2^{-n}|V(G_{n+1})|.\]

        Hence $(2)$ of Lemma~\ref{limit}
        holds for $\varepsilon=2^{-n}$, thus,
        the inverse limit will be a $d$-regular graphing without a non-zero antisymmetric circulation function and a measurable perfect matching.
	    Finally, since it has no non-zero antisymmetric circulation function it is essentially acyclic, that is, a treeing.
	\end{proof}
	
	\begin{proof}(of Theorem~\ref{coupling})
		The proof is similar to the previous one. We construct a regular sequence of bipartite graphs $\{G_n\}_{n=1}^{\infty}$, where $G_n$ is $d2^{n-1}$-regular. We start with $G_1=K_{d,d}$. In every step we construct the next graph from $G_n$ almost like in the previous proof, but in the end we replace every vertex by a couple of twin vertices and connect two new vertices if the original ones were adjacent. (In other words, we consider its tensor product with the graph on two adjacent vertices with a loop at both of them.) The assumptions of Lemma \ref{basic} but $(6)$ hold, like in the previous proof: the only difference is that the tensor product doubles the pre-image of every vertex, while the pre-image of an edge becomes four times bigger. But the relative sizes do not change. Hence the limit is a Markov space. Moreover, by $(7)$ of Lemma \ref{basic}, there exists a Borel bipartition $J=A \cup^* B$ such that $E(G) \subseteq A \times B \cup B \times A$ and $\eta(J \times J \setminus E(G))=0$. 

        We prove that the conditional probabilities are a.s. atomless. Now $(8)$ of Lemma \ref{basic} shows that the disintegration of $\eta$ is $\{ \eta_x \}_{x \in J}$ defined before Lemma \ref{basic}.
        Every measure $\eta_x$ is atomless: the measure of every element can be bounded by $\frac{1}{d2^{n-1}}$, the measure of one vertex in $G_n$, for every $n$.

        We define the potentials associated to the orientations of $G_n$ in the same way as in the previous proof,
        the value will be equal for twin vertices.
        Hence assumption $(1)$ of Lemma~\ref{limit} is trivial, while $(2)$ holds again for $\varepsilon=2^{-n}$: the tensor product doubles the vertex set and the set of bad vertices, too. So every circulation is symmetric. In particular, the only coupling satisfying $E(G)$ is $\eta$, so there is no such deterministic coupling by Lemma~\ref{limit}.
	\end{proof}
	
	\begin{proof}(of Theorem~\ref{action})
	    We may assume that $d$ is even: if $d$ was odd then a construction for $(d+1)$ would induce the required $(\mathbb{Z}/2\mathbb{Z})^{*d}$-action, and the Schreier graphing of the latter would be a subgraphing. We will use the graphing $G$ constructed in the proof of Theorem~\ref{main}. Following the notation of the proof the probability measure space where the group $(\mathbb{Z}/2\mathbb{Z})^{*d}$ acts will be the following subset of vertices:
		
        $$X=\{ x: \text{ }x_{n+1} \in V(G_{n+1}) \text{ is internal for every } n \}.$$

        Note that $X$ has positive measure, since $\frac{d^2|V(F_n)|}{|V(G_{n+1})|} \geq 1-2^{-n}$ holds for the proportion of internal vertices for every $n>1$, so $X$ is a probability measure space with the restriction of the node measure on $V(G)$.

        We will decompose a measurable subset of
        $E(G)$ into edge-disjoint paths such that the endvertices of every path are in $X$, and a.e. vertex of $X$ is the endvertex of exactly $d$ paths. And we will color these paths with $d$ colors such that different paths ending at the same vertex will have different color: such a coloring induces a.e. a $(\mathbb{Z}/2\mathbb{Z})^{*d}$-action on $X$ moving every vertex along the corresponding path, which is pmp, since the graph is pmp.
        We build this path decomposition by recursion using the following lemma.
        Define the sets $X_1=V(G_1)$ and $X_n=\{u: u \in V(G_n), u_i \text{ is an internal vertex of } G_i \text{ for }i=2, \dots n \}$ for $n \geq 2$. We will also allow closed paths, i.e., cycles, when partitioning the edge set of finite graphs. Cycles do not have endvertices.

        \begin{claim}
        For every $n$ the edge set $E(G_n)$ admits a partition into edge-disjoint paths and cycles such that every vertex of $X_n$ is the endvertex of $d$ paths, the endvertices of every path are in $X_n$ and there exists a $d$-coloring of the paths such that paths ending at the same vertex have different color. Moreover, every path in $G_n$ consisting of internal vertices only is the lift of a path in $G_{n-1}$.
        \end{claim}

        \begin{proof}
		We prove by induction on $n$. For $G_1=K_{d,d}$  every path consists of a single edge and we choose a proper $d$-coloring.
		When $G_n$ has been decomposed into paths and equipped with a $d$-coloring of these paths we try to lift this decomposition and coloring to every copy of $F_n$: two edges incident to the same vertex belong to the same path if their images under $f_n$ do, and the color of every path is the color of its image (more precisely, the color of the path containing its image). Note that the endvertices of the obtained paths are the pre-images of endvertices in $G_n$ (these are in $X_n$) or bad internal vertices.

        Next, we will concatenate the paths using edges incident to external vertices.
        Consider an external vertex $u \in V(G_{n+1})$. Concatenate every edge $(u,v)$ to the path ending at $v$ whose image is a subpath of the path containing $f_n(u,v)$. (If there is no such path then $v \in X_{n+1}$, and  the edge $(u,v)$ will become a path itself.) Match the paths ending at $u$ arbitrarily and merge them pairwise in order to get a new path decomposition, we use at this point that $d$ is even. Note that $f_n$ maps every neighbor of $u$ to the same vertex in $G_n$, hence the matched paths have the same color, and we give this color to the new merged paths. Observe that the endvertices in the path decomposition of $E(G_{n+1})$ are in $X_{n+1}$: for every internal endvertex of every path this either held throughout the process, or it was a bad vertex and the path was merged with an edge connecting this bad vertex to an external vertex. And there are no external endvertices left in the end.
        \end{proof}

		We will show that for a.e. $x \in X$ the paths ending at $x_n$ in $G_n$ have only internal vertices for all but finitely many $n$. If $x_n$ is contained by a path in $G_n$ for $x \in X$ then it should be the endvertex of the path.
		Every external vertex of $G_{n+1}$ is contained by exactly $d/2$ paths, every path has two endvertices
		and the proportion of external vertices is at most $2^{-n}$ for every $n$, hence, by the Borel-Cantelli lemma, for a.e. vertex $x \in X$ there are only  finitely many $n$ such that the vertex $x_n$ is contained by a path in $G_n$ with an external vertex.

        If a path in $G_{n+1}$ contains internal vertices only then it is the lift of a path in $G_n$.
        Hence for a.e. $x \in X$ there exists an $N$ such that for every $n>N$ the vertex $x_n \in G_n$ is the endvertex of exactly $d$ paths, and these are the lifts of paths in the decomposition of $G_{n-1}$.
        The colored decompositions of the finite graphs induce a colored decomposition of a measurable subset of $E(G)$ into colored paths, and a.e. $x \in X$
        will be the endvertex of $d$ paths.
		
		Given an antisymmetric bounded circulation function $f$ on the edge set of the Schreier graphing of the action of $(\mathbb{Z}/2\mathbb{Z})^{*d}$ on $X$ encoded by the colored path decomposition define the circulation function $g:E(G) \rightarrow \mathbb{R}$ as follows. For the vertices
		$x, y \in X$ and an edge $(x',y')$ on the path from $x$ to $y$ set $g(x',y')=f(x,y)$. Then $g$ is a bounded antisymmetric circulation function on $G$, hence by Theorem \ref{main} it is zero a.e., and
	so is $f$. Since every bounded antisymmetric
        circulation function is zero a.e., the action of $(\mathbb{Z}/2\mathbb{Z})^{*d}$ on $X$ is a.e. free.
	\end{proof}
	
	\section{Further directions}
	\label{s:further}

     \textbf{Perfect matchings}. It would be interesting to see which local properties of a graphing guarantee the existence of a measurable perfect matching. A group theoretic version of this problem is as follows.

    \begin{question}
    Which finitely generated groups admit a measurable perfect matching in their Schreier graphing of any pmp free ergodic action with respect to any finite symmetric set of generators?
    \end{question}

    This holds for Kazhdan groups (by \cite{CL}, and in the bipartite case by \cite{LN}, since Schreier graphs of their ergodic actions do always have a spectral gap) and for bipartite Schreier graphings of one-ended amenable groups \cite{BKS}. There are two-ended amenable counterexamples, see \cite{BHT}, \cite{Lmatching}. We do not know any non-amenable group without this property, though T\'oth and the author \cite{KT} have recently constructed a pmp free action of the group $(\mathbb{Z}/2\mathbb{Z})^{*d}$ which contains exactly $d$ measurable perfect matchings (with respect to the standard set of generators)! One wonders if the free group (with respect to the standard set of generators) is an example. The techniques of this paper fail to show this, since they only work for measurably bipartite graphings, and a measurably bipartite Schreier graphing of the free group clearly admits a measurable perfect matching.

    \textbf{Balanced orientations}.
    The question is similarly  intriguing for balanced orientations.

    \begin{question}
    Which finitely generated groups admit a measurable balanced orientation in the Schreier graphing of any pmp free ergodic action with respect to any finite symmetric set of generators?
    \end{question}


    We have seen that this
    holds for Kazhdan groups (by \cite{T}, since Schreier graphs of their ergodic actions do always have a spectral gap) and for one-ended amenable groups \cite{BKS}, though there are two-ended counterexamples \cite{C}, \cite{BHT}. Corollary \ref{balanced} gives a non-amenable counterexample, a  free pmp action of $(\mathbb{Z}/2\mathbb{Z})^{*d}$ (with respect to the standard set of generators).
    One can modify the construction to get an action that is also ergodic, but this would be beyond the goals of this paper.
    Note that the group $(\mathbb{Z}/2\mathbb{Z})^{*d}$ admits a set of generators consisting elements of infinite order (inducing a measurable balanced orientation for any pmp free action). These show that the answer could depend on the choice of generators, hence we ask if there exists a measurable balanced orientation with respect to any set of generators.
    While Schreier graphs of the free group always admit a measurable balanced orientation, T\'oth and the author have used the ideas of this paper to construct a free pmp action of the free group with one single measurable balanced orientation \cite{KT}.

    {\bf Separation for LCL problems.}
    Corollary \ref{balanced} gives for $d \geq 3$ a  free pmp action of the non-amenable group $(\mathbb{Z}/2\mathbb{Z})^{*d}$ without a measurable balanced orientation. On the other hand,
    Thornton's theorem shows that the Bernoulli shift of $(\mathbb{Z}/2\mathbb{Z})^{*d}$ admits a measurable balanced orientation. Hence the Bernoulli shift admits a
    finite measurable coloring whose colored local statistics can not be realized in a measurable way, though it can be approximated, since Bernoulli shifts are minimal (weakly contained) among free pmp actions by a theorem of Ab\'ert and Weiss \cite{AW}, see also Hatami, Lov\'asz and Szegedy \cite{HLS}. It would be interesting to see which LCL properties (see \cite{BCGGRV} for the definition) of the Bernoulli shift do not hold for every pmp free ergodic action of a group.

\end{document}